\documentclass[12pt]{smfart}

\usepackage[francais]{babel}
\usepackage{smfthm}
\usepackage{amssymb}
\usepackage{mathrsfs}

\newcommand{\val}{\mathrm{val}}
\newcommand{\Qp}{\mathbf{Q}_p}
\newcommand{\Zp}{\mathbf{Z}_p}
\newcommand{\Fp}{\mathbf{F}_p}
\newcommand{\ZZ}{\mathbf{Z}}
\newcommand{\OO}{\mathcal{O}}
\newcommand{\MM}{\mathfrak{m}}
\newcommand{\Fpbar}{\overline{\mathbf{F}}_p}
\newcommand{\Qpbar}{\overline{\mathbf{Q}}_p}
\newcommand{\dcris}{D_{\mathrm{cris}}}
\newcommand{\eps}{\varepsilon}
\newcommand{\ra}{\rightarrow}
\renewcommand{\phi}{\varphi}
\renewcommand{\projlim}{\varprojlim}
\newcommand{\Fil}{\mathrm{Fil}}
\renewcommand{\geq}{\geqslant}
\renewcommand{\leq}{\leqslant} 
\newcommand{\G}{\mathrm{GL}_2(\Qp)}
\newcommand{\B}{\mathrm{B}(\Qp)}
\newcommand{\K}{\mathrm{GL}_2(\Zp)}
\newcommand{\g}{\mathrm{Gal}(\Qpbar/\Qp)}
\newcommand{\calE}{\mathscr{E}}
\newcommand{\calO}{\mathscr{O}_{\mathscr{E}}}
\newcommand{\pr}{\mathrm{pr}}

\author[L. Berger]{Laurent Berger}
\address{C.N.R.S. \& I.H.\'E.S.\\
Le Bois-Marie\\
35 route de Chartres\\
91440 Bures-sur-Yvette \\ 
France}
\email{laurent.berger@ihes.fr}
\urladdr{www.ihes.fr/\~{}lberger/}

\author[C. Breuil]{Christophe Breuil}
\address{C.N.R.S. \& I.H.\'E.S.\\
Le Bois-Marie\\
35 route de Chartres\\
91440 Bures-sur-Yvette \\ 
France}
\email{breuil@ihes.fr}
\urladdr{www.ihes.fr/\~{}breuil/}

\title[Repr\'esentations cristallines en poids moyens]
{Sur la r\'eduction des repr\'esentations
  cristallines de dimension $2$ en poids moyens} 

\date{Juillet 2005}

\NumberTheoremsIn{subsection}

\begin{document}

\begin{abstract}
On calcule la r\'eduction modulo $p$
des repr\'esentations cristallines de dimension $2$ dont les poids de
Hodge-Tate sont $0$ et $k-1$ avec $k \in \{p+2,\cdots,2p-1\}$.
\end{abstract}

\begin{altabstract}
We compute the reduction modulo $p$ of $2$-dimensional crystalline
representations whose Hodge-Tate weights are $0$ and $k-1$ with $k \in
\{p+2,\cdots,2p-1\}$. 
\end{altabstract}

\maketitle

\setcounter{tocdepth}{2}

\tableofcontents

\setlength{\baselineskip}{18pt}

\section*{Introduction}
Soit $p$ un nombre premier et $L$ une
extension finie de $\Qp$, dont on note $\OO_L$, $\MM_L$ et $k_L$ 
l'anneau des entiers, l'id\'eal maximal, et le corps r\'esiduel. Si
$k$ est un entier $\geq 2$ et $a_p \in \MM_L$ on d\'efinit le
$\phi$-module filtr\'e $D_{k,a_p}$ par $D_{k,a_p} = L e \oplus L
f$ o\`u~: 
\[ \begin{cases} \phi(e) = p^{k-1} f \\
\phi(f) = -e + a_p f 
\end{cases}
\text{et}\quad
\Fil^i D_{k, a_p} = \begin{cases}
D_{k, a_p} & \text{si $i \leq 0$,} \\
L e & \text{si $1 \leq i \leq k-1$,} \\
0 & \text{si $i \geq k$.}
\end{cases} \]
Ce $\phi$-module filtr\'e est admissible, et  
on sait (par le th\'eor\`eme principal de \cite{CF})
qu'il existe alors une repr\'esentation cristalline
$V_{k,a_p}$ de $\g$
telle que $\dcris(V^*_{k,a_p}) = D_{k,a_p}$ (on passe au
dual pour que les notations soient compatibles avec celles de
\cite{Br2} et \cite{BLZ}; remarquons tout de m\^eme que 
l'on a $V_{k,a_p}^* =
V_{k,a_p}(1-k)$). 
Toute repr\'esentation cristalline 
absolument irr\'eductible de dimension $2$ de
$\g$ est la tordue par un caract\`ere cristallin d'une
$V_{k,a_p}$ avec $k \geq 2$.

Si $T_{k,a_p}$ est un $\OO_L$-r\'eseau de
$V_{k,a_p}$ stable par $\g$, alors la semi-simplifi\'ee de $k_L
\otimes_{\OO_L} 
T_{k,a_p}$ ne d\'epend pas du choix du r\'eseau et nous la notons
$\overline{V}_{k,a_p}$. La question se pose alors de donner une
formule pour $\overline{V}_{k,a_p}$ en termes de $k$ et $a_p$. Dans
\cite{Br2}, une formule conjecturale est donn\'ee pour 
$\overline{V}_{k,a_p}$ quand $2p \geq k \geq
2$. 

Quand $k \leq p$, cette conjecture suit imm\'ediatement de la 
{\og th\'eorie de Fontaine-Laffaille \fg} (cf. \cite{FL82}). 
Quand $k=p+1$ ou bien quand $2p-1 \geq
k \geq p+2$ et $\val(a_p) > 1$ 
(la valuation {\og $\val$ \fg} \'etant la valuation
$p$-adique) ou bien encore quand $k=2p$ et $\val(a_p)>2$, 
la conjecture est d\'emontr\'ee
dans \cite{BLZ}. L'objet de cet article est de d\'emontrer la
conjecture pour $2p-1 \geq k \geq p+2$ quand $0 < \val(a_p) \leq 1$ ce
qui en compl\`ete la d\'emonstration pour $k \leq 2p-1$. 
Notons $\omega$ le caract\`ere cyclotomique modulo $p$ et
$\mu_\lambda$ le caract\`ere non-ramifi\'e de $\g$ qui envoie
$\mathrm{Frob}_p^{-1}$ sur $\lambda$. On a alors
le th\'eor\`eme suivant qui vient compl\'eter \cite[proposition
6.2]{Br2} et \cite[theorem]{BLZ}~:

\begin{enonce*}{Th\'eor\`eme}
Pour $2p-1 \geq k \geq p+2$, 
la r\'eduction modulo $p$ des repr\'esentations $V_{k,a_p}$ est
donn\'ee par les formules ci-dessous. 
\begin{enumerate}
\item Pour $k=p+2$~:
\begin{enumerate}
\item si $1 > \val(a_p) > 0$, alors $\overline{V}_{k,a_p} =
  \mathrm{ind}(\omega_2^2)$. 
\item si $\val(a_p) = 1$, 
et si $\lambda$ est une racine du
polyn\^ome $\lambda^2 - \overline{a_p/p} \lambda +1 = 0$,
alors \[ \overline{V}_{k,a_p} = 
\begin{pmatrix} \omega\mu_{\lambda} & 0 \\ 
0 & \omega\mu_{\lambda^{-1}} \end{pmatrix}. \]
\end{enumerate}
\item Pour $2p-1 \geq k \geq p+3$~:
\begin{enumerate}
\item si $1 > \val(a_p) > 0$, alors $\overline{V}_{k,a_p} =
  \mathrm{ind}(\omega_2^{k-p})$. 
\item si $\val(a_p) = 1$, et
si $\lambda=\overline{a_p/p} \cdot (k-1)$, alors 
\[ \overline{V}_{k,a_p} = 
\begin{pmatrix} \omega^{k-2} \mu_{\lambda} & 0 \\ 
0 & \omega\mu_{\lambda^{-1}} \end{pmatrix}. \]
\end{enumerate}
\end{enumerate}
\end{enonce*}

Ce th\'eor\`eme est la r\'eunion des corollaires \ref{fin1},
\ref{fin2} et \ref{fin3}. La d\'emonstration est diff\'erente selon
que $\val(a_p)=1$ ou que $1 > \val(a_p) > 0$. Dans le premier cas, on
se borne \`a \'etendre les calculs de \cite{BLZ} qui traitaient du cas
$\val(a_p)>1$. Dans le deuxi\`eme cas, on utilise les id\'ees de
\cite{Br1,Br2} qui consistent \`a associer \`a $V_{k,a_p}$ une
repr\'esentation admissible de $\G$ et \`a v\'erifier, 
gr\^ace aux r\'esultats de \cite{CL} et de \cite{BB}, que cette
association (la {\og correspondance de Langlands $p$-adique continue
  \fg}) est compatible avec la r\'eduction modulo $p$ (dans tout cet
article, on fait l'abus de langage qui consiste \`a parler de {\og
  r\'eduction modulo $p$ \fg} quand on devrait 
plut\^ot parler de r\'eduction modulo $\MM_L$).

\section{Rappels et notations}

Comme l'objet de cet article est de calculer la r\'eduction modulo $p$
de certaines repr\'esen\-tations cristallines, nous supposons que le
lecteur est familier avec la notion de repr\'esenta\-tion cristalline et
de $\phi$-module filtr\'e. Nous faisons quelques rappels sur la
th\'eorie des $(\phi,\Gamma)$-modules, qui est essentielle pour la
suite, et sur quelques uns de ses prolongements. Pour des rappels
beaucoup plus d\'etaill\'es, nous renvoyons \`a \cite[\S4,5]{CL} 

\subsection{Repr\'esentations $p$-adiques et 
$(\varphi,\Gamma)$-modules}

Soit $\Gamma = \mathrm{Gal}(\Qp(\mu_{p^\infty})/\Qp)$ et $\eps :
\Gamma \ra \Zp^\times$ le caract\`ere cyclotomique et $\calO$ l'anneau
$\calO = \{ \sum_{i \in \ZZ} a_i X^i \}$ o\`u $a_i \in \OO_L$
et $a_{-i} \ra 0$ quand $i \ra \infty$. On munit cet anneau d'un
frobenius $\OO_L$-lin\'eaire $\phi$ d\'efini par $\phi(X)=(1+X)^p-1$
et d'une action $\OO_L$-lin\'eaire de $\Gamma$ donn\'ee par $\gamma(X)
= (1+X)^{\eps(\gamma)}-1$ si $\gamma \in \Gamma$. Un
$(\phi,\Gamma)$-module \'etale est un $\calO$-module $D$ de type fini
muni d'un frobenius semi-lin\'eaire $\phi$ et d'une action de $\Gamma$
semi-lin\'eaire continue 
et commutant \`a $\phi$. Rappelons que Fontaine a
construit dans \cite[A.3.4]{F90}
un foncteur $T \mapsto D(T)$ qui \`a toute
$\OO_L$-repr\'esentation de $\g$ associe un $(\phi,\Gamma)$-module
\'etale et que ce foncteur est une \'equivalence de cat\'egories. 

\subsection{Modules de Wach}

Un module de Wach est un $\OO_L[\![X]\!]$-module $N$ muni d'un frobenius
$\phi$ et d'une action de $\Gamma$ semi-lin\'eaire continue 
et commutant \`a $\phi$, satisfaisant les conditions
suivantes~: 
\begin{enumerate}
\item $N$ est libre de rang fini sur $\OO_L[\![X]\!]$;
\item le groupe $\Gamma$ agit trivialement sur $N/X$;
\item il existe $h \geq 0$ tel que $N/\phi^*(N)$ est tu\'e par $q^h$
  o\`u $q=\phi(X)/X$. 
\end{enumerate}
Le plus petit entier $h$ 
v\'erifiant (3) ci-dessus
est appel\'e la hauteur du module de Wach $N$.

Si $T$ est une $\OO_L$-repr\'esentation sans torsion telle que 
$V = L \otimes_{\OO_L} T$ est cristalline \`a poids de Hodge-Tate $\leq 0$,
alors l'un des principaux r\'esultats de \cite{Be1} est  
qu'il existe un (unique) module de Wach $N(T)$ tel que $D(T) = \calO
\otimes_{\OO_L[\![X]\!]} N(T)$. Posons alors $N(V) = L \otimes_{\OO_L} N(T)$
et $\Fil^i N(V) = \{ x \in N(V)$ tels que $\phi(x) \in q^i N(V)\}$. Le
$L$-espace vectoriel $N(V)/X$ h\'erite de la filtration induite et du
frobenius $\phi$, ce qui en fait un $\phi$-module filtr\'e. Le
th\'eor\`eme III.4.4 de \cite{Be1} nous dit alors que $N(V)/X \simeq
\dcris(V)$. On en d\'eduit facilement le fait suivant~: 
\begin{lemm}\label{ggg}
Si $N$ est un module
de Wach et si $V$ est une repr\'esentation cristalline 
telle que $\dcris(V)=N/X$, alors le $(\phi,\Gamma)$-module 
$D(V)$ associ\'e \`a $V$ est isomorphe \`a $\calE \otimes_{\OO_L[\![X]\!]} N$.
\end{lemm}

\subsection{L'op\'erateur $\psi$ et le module
  $D^\sharp(V)$}\label{rap3} 

L'anneau $\calO$ est un $\phi(\calO)$-module libre de rang $p$,
dont une base est donn\'ee par $\{(1+X)^i\}_{0 \leq i \leq p-1}$. Si
$x \in \calO$, on peut donc \'ecrire $x=\sum_{i=0}^{p-1} (1+X)^i
\phi(x_i)$ et on d\'efinit un op\'erateur $\psi : \calO \rightarrow
\calO$ par la formule $\psi(x) = x_0$ si $x=\sum_{i=0}^{p-1} (1+X)^i
\phi(x_i)$. 

Si $D$ est un $(\phi,\Gamma)$-module \'etale sur $\calO$, alors Colmez
a d\'efini dans \cite[\S 4.5]{CL} un sous-$\OO_L[\![X]\!]$-module
$D^\sharp$ de $D$. Si $D = \calO \otimes_{\OO_L[\![X]\!]} N$ o\`u $N$
est un module de Wach de hauteur $h$, 
alors $D^\sharp$ est caract\'eris\'e 
par les propri\'et\'es suivantes 
(cf. \cite[\S 4]{CL} et \cite[\S 3]{BB})~: 
\begin{enumerate}
\item $D^\sharp \subset X^{-h-1}N$;
\item  quels que soient $x \in D$ et $j \geq 0$, il existe $n(x,j)
  \geq 0$ tel que $\psi^n(x) \in D^{\sharp} + p^j D$ si $n \geq n(x,j)$;
\item l'op\'erateur $\psi$ induit une surjection de $D^{\sharp}$ sur
  lui-m\^eme.
\end{enumerate}
Le $\OO_L[\![X]\!]$-module $D^\sharp$ est donc stable par $\psi$ et
aussi sous l'action de $\Gamma$. Nous l'utiliserons un peu plus loin,
au paragraphe \ref{defllpc}.

\section{Calcul de la r\'eduction~: le cas $\mathrm{val}(a_p) = 1$}

Dans ce chapitre, on construit les modules de Wach associ\'es aux
repr\'esentations $V_{k,a_p}^*$ pour $2p-1 \geq k \geq p+2$ et
$\val(a_p)=1$ puis on utilise les formules explicites ainsi obtenues
pour calculer $\overline{V}_{k,a_p}$.  

\subsection{Construction du module de Wach}

La construction des modules de Wach associ\'es aux repr\'esentations
$V_{k,a_p}^*$ est la m\^eme que celle que l'on a donn\'ee dans
\cite{BLZ} (mais attention au fait que les notations sont
l\'eg\`erement diff\'erentes). 
Nous en rappelons ici les points essentiels. Pour $n \geq
1$, on pose $q_n = \phi^{n-1}(\phi(X)/X)$ ce qui fait en particulier
que $q_1=q=\phi(X)/X$ et on d\'efinit deux s\'eries $\lambda_+$ et
$\lambda_-$ par les formules~:
\[ \lambda_+ = \prod_{n \geq 0} \frac{\phi^{2n+1}(q)}{p} = 
\frac{q_2}{p} \times \frac{q_4}{p} \times \frac{q_6}{p} \times \cdots 
\qquad\text{et}\qquad
\lambda_-  = \prod_{n \geq 0} \frac{\phi^{2n}(q)}{p} =
\frac{q_1}{p} \times \frac{q_3}{p} \times \frac{q_5}{p} \times \cdots \]

Puisque l'on suppose que $\val(a_p)=1$, 
la proposition suivante r\'esulte du (4) de \cite[proposition
3.1.1]{BLZ}~:

\begin{prop}\label{lpm}
Si l'on \'ecrit $a_p (\lambda_-/\lambda_+)^{k-1} = \sum_{i \geq
    0} \alpha_i X^i$, alors $\val(\alpha_i) \geq 0$ pour
  $i\in\{0,\cdots,k-2\}$. 
\end{prop}

On pose alors $\alpha=\alpha_0+\alpha_1 X + \cdots + \alpha_{k-2}
X^{k-2}$ ainsi que $g_\pm=\lambda_\pm/\gamma(\lambda_\pm)$
et on d\'efinit une matrice $P \in 
\mathrm{M}_2(\OO_L[\![X]\!])$ et, pour tout $\gamma \in
\Gamma$, une matrice $G_\gamma^{(k-1)} \in \mathrm{Id} + X \cdot
\mathrm{M}_2(\OO_L[\![X]\!])$ par les formules~: 
\[ P = 
\begin{pmatrix}
0 & -1 \\
q^{k-1} & \alpha
\end{pmatrix}
\quad\text{et}\quad
G^{(k-1)}_\gamma = 
\begin{pmatrix}
g_+^{k-1} & 0 \\
0 & g_-^{k-1}
\end{pmatrix}. \]

\begin{prop}\label{unikgam}
Si $\gamma \in \Gamma$, alors il existe 
une unique matrice $G_\gamma \in 
\mathrm{M}_2(\OO_L[\![X]\!])$ telle que~: 
\begin{enumerate}
\item $P \phi(G_\gamma) = G_\gamma \gamma(P)$; 
\item $G_\gamma \equiv
G_\gamma^{(k-1)} \mod{X^{k-1} \cdot \mathrm{M}_2(\OO_L[\![X]\!])}$.
\end{enumerate}
\end{prop}

\begin{proof}
Montrons tout d'abord l'unicit\'e de la matrice $G_\gamma$. Si 
$G_\gamma$ et $G'_\gamma$ sont deux matrices 
satisfaisant les conditions de la
proposition et si l'on pose $H=G'_\gamma G^{-1}_\gamma$, alors un
petit calcul montre que $H P = P \phi(H)$ ce qui fait que si l'on
\'ecrit $H = \mathrm{Id} + H_\ell X^\ell + \cdots$ avec $H_\ell \in
\mathrm{M}_2(\OO_L)$ et $H_\ell \neq 0$ et que $P_0$ d\'enote le
coeffcient constant de $P$, alors on a $H_\ell P_0  = p^\ell P_0 H_\ell$
ce qui implique que $P_0$ a deux valeurs propres dont le quotient est
$p^\ell$. Comme $\ell \geq k-1$ par la condition (2) 
de la proposition, ce n'est pas possible et 
$H=\mathrm{Id}$ et donc $G_\gamma = G_\gamma'$. 

La d\'emonstration de l'existence de $G_\gamma$ est semblable \`a
celle qui est donn\'ee dans la preuve de \cite[proposition
3.1.3]{BLZ}, nous en rappelons ici les points essentiels. Un calcul
direct montre qu'il existe $R^{(k-1)} \in
\mathrm{M}_2(\OO_L[\![X]\!])$ telle que~:
\[ G^{(k-1)}_\gamma - P \phi(G^{(k-1)}_\gamma)
\gamma(P^{-1}) = X^{k-1} R^{(k-1)}. \]
La fin de la 
d\'emonstration consiste \`a montrer par r\'ecurrence sur $\ell \geq
k$ qu'il existe deux matrices $R^{(\ell)}$ et $G_\gamma^{(\ell)}$
\`a coefficients dans $\OO_L[\![X]\!]$ telles que~: 
\begin{enumerate}
\item $G^{(\ell)}_\gamma \equiv G^{(\ell-1)}_\gamma \mod{X^{\ell-1}}$;
\item  $G^{(\ell)}_\gamma - P \phi(G^{(\ell)}_\gamma)
\gamma(P^{-1}) = X^{\ell} R^{(\ell)}$.
\end{enumerate}
Ceci se d\'emontre par approximations successives, et il suffit alors
de prendre $G_\gamma = \lim_{\ell \ra + \infty} G_\gamma^{(\ell)}$. 
\end{proof}

On d\'efinit alors un module de Wach $N_{k,a_p} 
= \OO_L[\![X]\!] e \oplus \OO_L[\![X]\!] f$ en d\'ecidant que les
matrices de $\phi$ et de $\gamma \in \Gamma$ dans la base $\{e,f\}$
sont donn\'ees par $P$ et $G_\gamma$. Le (2) de la proposition
\ref{unikgam} montre que $G_{\gamma \eta} = G_\gamma \gamma(G_\eta)$
et donc que cela d\'efinit bien une action semi-lin\'eaire de
$\Gamma$, et le (1) de la proposition implique que
$\phi$ commute \`a cette action du groupe $\Gamma$.

\begin{prop}\label{verifwach}
Le $(\phi,\Gamma)$-module $\calE \otimes_{\OO_L[\![X]\!]} N_{k,a_p}$ est isomorphe
\`a $D(V_{k,a_p}^*)$. 
\end{prop}

\begin{proof}
Par le lemme \ref{ggg}, il suffit de v\'erifier que le
$\phi$-module filtr\'e $N_{k,a_p} / X$ est bien isomorphe \`a
$D_{k,a_p}$ ce que nous laissons en exercice au lecteur (c'est 
\cite[proposition 3.2.4]{BLZ}). 
\end{proof}

\subsection{Calcul de la r\'eduction}

L'objet de ce paragraphe est d'utiliser les formules explicites du
paragraphe pr\'ec\'edent pour calculer la r\'eduction modulo $p$ des
repr\'esentations $V_{k,a_p}$ quand $2p-1 \geq k \geq p+2$ et
$\val(a_p)=1$.  

\begin{lemm}\label{alphamod}
Si $\beta=\overline{a_p/p} \cdot (k-1)$, alors 
il existe $u(X) \in 1 + X \cdot k_L[\![X]\!]$ tel que
$\overline{\alpha} = \beta \cdot u(X) \cdot X^{p-1}$.
\end{lemm}

\begin{proof}
Rappelons que $\alpha$ a \'et\'e d\'efini comme la 
partie de degr\'e $\leq k-2$ de la s\'erie~: 
\[ a_p \cdot
\left( \frac{(q_1/p) \cdot (q_3/p) \cdot (q_5/p) \cdots }
{(q_2/p) \cdot (q_4/p) \cdots} \right)^{k-1} = 
 a_p \cdot \left( \frac{q_1}{p} \right)^{k-1}
\left(  \frac{ (q_3/p) \cdot (q_5/p) \cdots }
{(q_2/p) \cdot (q_4/p) \cdots} \right)^{k-1}. \] 
On sait que $q_1(0)=p$ et que $q_1 \equiv X^{p-1} \mod{p}$ 
et donc que si $n \geq 1$, alors $q_n(0)=p$ et 
$q_n \equiv X^{p^{n-1}(p-1)} \mod{p}$ ce qui fait que la 
partie de degr\'e $\leq k-2$ de la s\'erie $(q_n/p)^{\pm 1}$ est \`a
coefficients dans $\Zp$ si $n \geq 2$. On en conclut
qu'il existe une s\'erie $v(X) \in 1 + X \cdot \Zp[\![X]\!]$ 
telle que $\alpha$ est la partie de degr\'e $\leq k-2$ de la s\'erie~:
\[ a_p \left(1+\frac{X^{p-1}}{p}\right)^{k-1} v(X)  = a_p \left( 1 + (k-1)
 \frac{X^{p-1}}{p} \right)  v(X) + O(X^{2(p-1)}). \]
Le lemme en r\'esulte imm\'ediatement, puisque $k-2 < 2(p-1)$.
\end{proof}

Nous allons maintenant passer au calcul proprement dit de la
r\'eduction modulo $p$. 

\subsubsection*{Le cas $2p-1 \geq k \geq p+3$}
Le cas $k=p+2$ se comporte de mani\`ere l\'eg\`erement diff\'erente du
cas $k \geq p+3$, et nous commen\c{c}ons par traiter ce dernier. Nous
allons montrer que le $(\phi,\Gamma)$-module $D(\overline{T}^*_{k,a_p})$ 
contient un sous-objet de rang $1$ que nous identifions. Dans ce
paragraphe, on pose $\lambda=\beta=\overline{a_p/p} \cdot (k-1)$. 

\begin{lemm}\label{solz}
Si $k \geq p+3$, alors il existe une unique s\'erie 
$z \in 1 + X \cdot k_L[\![X]\!]$ telle que~:
\[ z - u \phi(z) + \frac{X^{k-p-2}}{\lambda^2} \phi^2(z) = 0 \]
\end{lemm}

\begin{proof}
Un calcul imm\'ediat montre 
d'une part que l'application $z \mapsto u \phi(z) - (X^{k-p-2} / \lambda^2)
\phi^2(z)$ pr\'eserve $1 + X \cdot k_L[\![X]\!]$ (c'est l\`a qu'on 
utilise le fait que $k-p-2 \geq 1$)
et d'autre part que si
$r \geq 1$ est tel que $X^r$ divise $z$, alors $X^{pr}$ divise $u
\phi(z) - (X^{k-p-2} / \lambda^2) \phi^2(z)$. Ceci implique que
l'application $z \mapsto u \phi(z) - (X^{k-p-2} / \lambda^2)
\phi^2(z)$ est contractante 
pour la topologie $X$-adique
sur $1 + X \cdot k_L[\![X]\!]$ et donc qu'elle y
admet un unique point fixe.
\end{proof}

On pose alors~: \[ \delta = - \frac{\phi(z)}{\lambda}
\frac{e}{X^p} + z
\frac{f}{X} \in D(\overline{T}^*_{k,a_p}) = k_L(\!(X)\!)
\otimes_{k_L[\![X]\!]} N_{k,a_p}. \]

\begin{prop}\label{souze}
La $k_L(\!(X)\!)$-droite engendr\'ee par $\delta$ est 
un sous-$(\phi,\Gamma)$-module de $D(\overline{T}^*_{k,a_p})$
qui correspond \`a la sous-repr\'esentation
$\omega^{-1}\mu_{\lambda}  \subset \overline{T}^*_{k,a_p}$. 
\end{prop}

\begin{proof}
Un calcul imm\'ediat montre que $\phi(\delta)=\lambda \delta$. 
Si $\gamma \in \Gamma$, alors rappelons que 
$\gamma(e) - g_+^{k-1} e$ et $\gamma(f)
- g_-^{k-1} f$ appartiennent \`a $X^{k-1} \cdot 
k_L[\![X]\!] e \oplus X^{k-1}
\cdot k_L[\![X]\!] f$ et donc que si $x,y \in k_L[\![X]\!]$, alors~:
\[ \gamma\left(x \frac{e}{X^p} + y \frac{f}{X}
\right) = x' \frac{e}{X^p} + y' \frac{f}{X}, \]
avec \[ x' \equiv \gamma(x) \frac{X^p}{\gamma(X^p)} g_+^{k-1} 
\equiv x \omega(\gamma)^{-1} \mod{X} \]
et \[ y' \equiv \gamma(y) \frac{X}{\gamma(X)} g_-^{k-1} 
\equiv y \omega(\gamma)^{-1} \mod{X}. \]
Si l'on pose $x=-\phi(z)/\lambda$ et $y=z$, 
et que $x'$ et $y'$ sont d\'efinis comme ci-dessus,
alors~: 
\[ \phi \left(x' \frac{e}{X^p} + y' \frac{f}{X}
\right) = \phi \circ \gamma (\delta) = \gamma \circ \phi(\delta) 
= \gamma( \lambda \delta) = \lambda 
\left(x' \frac{e}{X^p} + y' \frac{f}{X} \right)
\] ce qui fait que
$\omega(\gamma) y'$ satisfait les conditions du
lemme \ref{solz} et donc que $y' = \omega(\gamma)^{-1} z$ et que $x' = 
-\omega(\gamma)^{-1} \phi(z)/\lambda$ ce qui fait que
$\gamma(\delta)=\omega(\gamma)^{-1} \delta$.

La droite $k_L(\!(X)\!) \delta$ est donc un $(\phi,\Gamma)$-module de rang $1$
avec $\phi(\delta)=\lambda \delta$ et
$\gamma(\delta)=\omega(\gamma)^{-1} \delta$; ce
$(\phi,\Gamma)$-module correspond au caract\`ere $\omega^{-1}
\mu_{\lambda}$.  
\end{proof}

La proposition \ref{souze} montre que $\overline{T}^*_{k,a_p}$
contient comme sous-repr\'esentation le caract\`ere 
$\omega^{-1}\mu_{\lambda}$ et donc que $\overline{V}_{k,a_p}$ 
(\'etant semi-simple) contient comme sous-repr\'esentation le caract\`ere 
$\omega\mu_{\lambda^{-1}}$. Comme $\det  \overline{V}_{k,a_p} = 
\omega^{k-1}$, on en d\'eduit finalement que~: 
\begin{coro}\label{fin1}
Si $2p-1 \geq k \geq p+3$ et $\val(a_p) = 1$ et
$\lambda=\overline{a_p/p} \cdot (k-1)$, alors~:
\[ \overline{V}_{k,a_p} = 
\begin{pmatrix} \omega^{k-2} \mu_{\lambda} & 0 \\ 
0 & \omega\mu_{\lambda^{-1}} \end{pmatrix}. \]
\end{coro}

\subsubsection*{Classe de l'extension pour $k=p+3$ et $\lambda=\pm 1$}

On conserve les notations du paragraphe pr\'ec\'edent; quand $k=p+3$ et $\lambda=\pm 1$, on a $\overline{V}_{k,a_p} = (\omega \oplus 1) \otimes \omega^{-2} \mu_{\lambda}$ et la question se pose alors de savoir si, avant semi-simplification, la r\'eduction modulo $p$ de $V_{k,a_p}$ est peu ramif\'ee ou tr\`es ramifi\'ee. La r\'eponse \`a cette question est donn\'ee par le th\'eor\`eme suivant~:

\begin{theo}\label{pourkisin}
Si $k=p+3$ et $\lambda=\pm 1$, alors~:
\[ \overline{T}^*_{k,a_p} = 
\begin{pmatrix} \omega  & \star \\ 
0 & 1 \end{pmatrix} \otimes \omega^{-2} \mu_{\lambda}, \]
o\`u $\star$ est non-trivial et peu ramifi\'e. 
\end{theo}

\begin{proof}
Nous ne donnons ici que les grandes lignes de la d\'emonstration de ce th\'eor\`eme.
Rappelons que si l'on pose comme ci-dessus~: \[  \delta = - \frac{\phi(z)}{\lambda}
\frac{e}{X^p} + z
\frac{f}{X}, \]
alors $\phi(\delta)  = \lambda \delta$ et $\gamma(\delta)  = \omega(\gamma)^{-1} \delta$.  Un calcul montre que~:
\[ \phi \left( \frac{e}{z X^{p+1}} \right) = \frac{1}{\phi(z)z X} \left(
z \frac{f}{X} -  \frac{\phi(z)}{\lambda}
\frac{e}{X^p} \right) + \frac{1}{\lambda} \frac{e}{z X^{p+1}}, \]
et donc que la matrice de $\phi$ dans la base $(\delta,e/z X^{p+1})$ est  donn\'ee par~: \[ \mathrm{Mat}(\phi) =  \begin{pmatrix}  \lambda & \frac{1}{z \phi(z)} \frac{1}{X} \\
0 & \lambda \end{pmatrix}. \]
De m\^eme, on v\'erifie que la matrice de $\gamma$ est donn\'ee par~:
\[ \mathrm{Mat}(\gamma) =  \begin{pmatrix}  \omega(\gamma)^{-1} & X^2 v_\gamma(X) \\ 0 & \omega(\gamma)^{-2} \end{pmatrix}, \]
o\`u $v_\gamma(X) \in k_L[\![X]\!]$. 

Apr\`es torsion par $\omega^2 \mu_{\lambda}$, on obtient donc une extension du $(\phi,\Gamma)$-module trivial par celui de $\Fp(1)$, extension donn\'ee par la classe~: \[ \mathrm{cl} \left( \frac{1}{X} (\lambda+\cdots), X^2 (1+\cdots)\right)   \in H^1\left(C_{\phi,\gamma}(\Qp,\Fp(1))\right), \] dans les notations de \cite[\S I.4]{CC}.
Cette classe est non-triviale par les lemmes I.5.2 et I.5.5 de \cite{CC}
car l'image de $-\psi(X^{-1}(\lambda+\cdots))$ est non-nulle dans $D(\Fp(1))/(\psi-1)$.
Pour terminer, il faut faire le calcul explicite de l'image de l'application de Kummer en termes de $(\phi,\Gamma)$-modules, ce qui est fait en partie dans \cite[\S V.3]{CC} et en partie dans \cite[\S 2.1]{Bn1}; on en conclut notamment que les extensions tr\`es ramifi\'ees correspondent aux classes $\mathrm{cl}(x,y)$ telles que la s\'erie formelle $y$ a un r\'esidu non-nul en $X=0$, ce qui n'est pas le cas ici.
\end{proof}

\subsubsection*{Le cas $k=p+2$}
Revenons au calcul de $\overline{V}_{k,a_p}$; 
dans le cas $k=p+2$, l'analogue du 
lemme \ref{solz} du paragraphe pr\'ec\'edent est faux
(parce que $X^{k-p-2}=1$) et il faut proc\'eder un petit 
peu diff\'erement. 

Remarquons que le $k_L[\![X]\!]$-module engendr\'e par
$e/X^p$ et $f/X$ est stable par
$\phi$ et que la matrice $Q$ de $\phi$
dans cette base v\'erifie $Q \in \mathrm{GL}_2(k_L[\![X]\!])$ 
(ce qui n'\'etait pas le cas si $2p-1 \geq k \geq p+3$). 

\begin{lemm}\label{dwork}
Si $Q \in \mathrm{GL}_2(k_L[\![X]\!])$, alors il existe $M \in
\mathrm{Id} + X \cdot
\mathrm{M}_2(k_L[\![X]\!])$ telle que $M^{-1} Q \phi(M) = Q(0)$.
\end{lemm}

\begin{proof}
Si on \'ecrit $Q=\sum_{i \geq 0} Q_i X^i$ et $M=\sum_{i \geq 0} M_i X^i$
avec $M_0=\mathrm{Id}$, alors on v\'erifie ais\'ement 
que les $M_i$ sont donn\'es par
r\'ecurrence par la formule 
$M_i  =   ( \sum_{j=0}^{\lfloor i/p \rfloor} Q_{i-pj} M_j
) Q_0^{-1}$.
\end{proof}

Un calcul facile montre que~:
\[ Q(0) = \begin{pmatrix} 0 & -1 \\ 1 & \beta \end{pmatrix}. \]
Si $\lambda$ est une racine dans $k_L$ du polyn\^ome $\lambda^2 -
\beta \lambda + 1 = 0$ (on suppose que $k_L$ contient cette racine, ce
qui est possible quitte \`a agrandir $L$), alors  
le lemme \ref{dwork} ci-dessus montre que le 
$k_L[\![X]\!]$-module engendr\'e par $e/X^p$ et
$f/X$ a une base $e',f'$ dans laquelle on a
$\phi(e')=\lambda e'$ et $\phi(f')= \lambda^{-1} f'$. 

Le fait que si $\gamma \in \Gamma$, alors
$\gamma(e) - g_+^{k-1} e$ et $\gamma(f)
- g_-^{k-1} f$ appartiennent \`a $X^{k-1} \cdot k_L[\![X]\!] e \oplus X^{k-1}
\cdot k_L[\![X]\!] f$ implique par ailleurs que la matrice de $\gamma$ dans la
base $e',f'$ est scalaire et diagonale. Comme
$\gamma(e/X^p) - \omega(\gamma)^{-1} e/X^p$ et $\gamma(f/X) -
\omega(\gamma)^{-1} f/X$ appartiennent \`a $X \cdot k_L[\![X]\!] (e/X^p) \oplus
X \cdot k_L[\![X]\!] (f/X)$, on voit que $\gamma$ agit par
$\omega(\gamma)^{-1}$, et donc finalement que $\overline{T}^*_{k,a_p}
=  \omega^{-1} \mu_\lambda \oplus \omega^{-1}
\mu_{\lambda^{-1}}$ et donc que~: 
\begin{coro}\label{fin2}
Si $k=p+2$ et $\val(a_p) = 1$ et $\lambda$ est une racine du
polyn\^ome $\lambda^2 - \overline{a_p/p} \lambda +1 = 0$, alors~:
\[ \overline{V}_{k,a_p} = 
\begin{pmatrix} \omega\mu_{\lambda} & 0 \\ 
0 & \omega\mu_{\lambda^{-1}} \end{pmatrix}. \]
\end{coro}

\section{Calcul de la r\'eduction~: le cas  $0 < \mathrm{val}(a_p)  < 1$}

Le calcul de $\overline{V}_{k,a_p}$ quand $0 < \val(a_p) < 1$
se fait par des m\'ethodes compl\`etement diff\'erentes de celles du 
paragraphe pr\'ec\'edent. On utilise ici la {\og correspondance de 
Langlands $p$-adique continue \fg}. 

\subsection{Une repr\'esentation du Borel sup\'erieur}\label{defllpc} 
Soient $N_{k,a_p}$ le module de Wach du dual $T^*_{k,a_p}$ 
d'un r\'eseau $T_{k,a_p}$ de $V_{k,a_p}$,
$D^\sharp(T_{k,a_p})$ le sous-$\OO_L[\![X]\!]$-module
de $X^{-k} N_{k,a_p}$ d\'efini au paragraphe \ref{rap3} et
$\projlim_\psi D^\sharp(T_{k,a_p})$ l'ensemble des suites
$\{v_n\}_{n \geq 0}$ d'\'el\'ements $v_n \in D^\sharp(T_{k,a_p})$
telles que $\psi(v_{n+1})=v_n$. 
On \'ecrit $\projlim_\psi D^\sharp(V_{k,a_p})$ pour $L \otimes_{\OO_L}
\projlim_\psi D^\sharp(T_{k,a_p})$.  
On d\'efinit une action du 
sous-groupe de Borel sup\'erieur $\B$ de $\G$~:
\[ \B = \left\{ \begin{pmatrix} * & * \\ 0 & * \end{pmatrix} \right\}
\subset \G \]
sur $\projlim_\psi D^\sharp(V_{k,a_p})$ de la mani\`ere suivante. 
Tout \'el\'ement $g \in \B$ peut s'\'ecrire comme produit~:
\[ g = 
\begin{pmatrix} x & 0 \\ 0 & x \end{pmatrix} \cdot
\begin{pmatrix} 1 &  0 \\ 0 & p^j  \end{pmatrix} \cdot 
\begin{pmatrix} 1 &  0 \\ 0 & a  \end{pmatrix} \cdot 
\begin{pmatrix} 1 &  z \\ 0 & 1 \end{pmatrix}, \]
o\`u $x \in \Qp^\times$, $j\in\ZZ$, $a \in \Zp^{\times}$ et $z
\in \Qp$. Si $v=\{v_n\}_{n\geq 0} \in \varprojlim_{\psi}
D^{\sharp}(V_{k,a_p})$, alors on pose pour $n\geq 0$~:
\begin{align*}
\left( \begin{pmatrix} x & 0 \\ 0 & x \end{pmatrix} \cdot v \right)_n
& = x_0^{k-2} v_n \text{ o\`u $x=p^{\val(x)} x_0$;} \\
\left( \begin{pmatrix} 1 &  0 \\ 0 & p^j  
\end{pmatrix} \cdot v \right)_n & = v_{n-j} = \psi^j(v_n); \\
\left( \begin{pmatrix} 1 &  0 \\ 0 & a  
\end{pmatrix} \cdot v \right)_n & = \gamma_a(v_n) 
\text{ o\`u $\gamma_a \in \Gamma$ est tel que $\eps(\gamma_a) = a$;}\\
\left( \begin{pmatrix} 1 &  z \\ 0 & 1  
\end{pmatrix} \cdot v \right)_n & = 
\psi^m((1+X)^{p^{n+m} z} v_{n+m}),\ n+m\geq -{\rm val}(z).
\end{align*}

Si $\chi$ est un caract\`ere cristallin de $\g$, on note
$V_{k,a_p,\chi} = V_{k,a_p} \otimes \chi$. 
Dans \cite{Br2}, il est construit 
des repr\'esentations localement alg\'ebriques de $\G$ sur $L$ 
not\'ees $\Pi_{k,a_p,\chi}$ ainsi que des $\OO_L$-r\'eseaux
$\Theta_{k,a_p,\chi} \subset \Pi_{k,a_p,\chi}$ stables par $\G$ pour $k
\leq 2p$. Soit 
$\Pi(V_{k,a_p,\chi})$ le compl\'et\'e $p$-adique de $\Pi_{k,a_p,\chi}$
par rapport \`a $\Theta_{k,a_p,\chi}$. La notation
$\Pi(V_{k,a_p,\chi})$ est justifi\'ee par le fait que dans \cite{BB}
il est d\'emontr\'e,  
en utilisant l'id\'ee principale de \cite{CL}, 
que, lorsque le Frobenius sur $D_{k,a_p}$ est semi-simple, alors on a un isomorphisme topologique $\B$-\'equivariant
entre le Banach dual $\Pi(V_{k,a_p,\chi})^*$ (muni de sa topologie
faible) et $\projlim_\psi D^\sharp(V_{k,a_p,\chi})$. 

D'autre part, les r\'eductions
$k_L \otimes_{\OO_L} \Theta_{k,a_p,\chi}$ ont \'et\'e 
d\'etermin\'ees explicitement dans \cite{Br2}~: 
dans les cas qui nous int\'eressent ($p+2 \leq k \leq 2p$ et $0
< \mathrm{val}(a_p)  < 1$) ce sont des repr\'esentations
supersinguli\`eres de $\G$. Ce sont ces formules explicites et le lien
entre $\Pi(V_{k,a_p,\chi})$ et  $\projlim_\psi
D^\sharp(V_{k,a_p,\chi})$ qui vont nous permettre de calculer
$\overline{V}_{k,a_p}$. 

\subsection{Repr\'esentations modulaires et supersinguli\`eres}
Rappelons que dans \cite{Br1}, on donne la liste de ces
repr\'esentations supersinguli\`eres de $\G$, que nous allons rappeler
pour la convenance du lecteur. Si $r \in \{
0,\cdots,p-1\}$ et si $\chi : \Qp^\times \ra k_L^{\times}$ est un
caract\`ere continu, on pose~: 
\[ \pi(r,0,\chi) = \left[ \left( \mathrm{ind}_{\K\Qp^\times}^{\G}
    \mathrm{Sym}^r k_L^2 \right) / T \right] \otimes (\chi \circ
\det), \] 
o\`u $T$ est un certain op\'erateur de Hecke. 
Par \cite[th\'eor\`eme 1.3]{Br1}, les entrelacements entre les
$\pi(r,0,\chi)$ sont les suivants~: 
\begin{align*} 
\pi(r,0,\chi) \simeq & \pi(r,0,\chi \mu_{-1}) \\
\pi(r,0,\chi) \simeq & \pi(p-1-r,0,\chi \omega^r) \\
\pi(r,0,\chi) \simeq & \pi(p-1-r,0,\chi \omega^r \mu_{-1}).
\end{align*} 

On peut d'autre part
faire une liste des repr\'esentations 
absolument irr\'eductibles de $\g$ de 
dimension $2$ sur $k_L$. Si $r \in \{
0,\cdots,p-1\}$ et si $\chi : \Qp^\times \ra k_L^\times$ est un
caract\`ere continu, que l'on identifie \`a un caract\`ere continu de
$\g$ via le corps de classes (normalis\'e pour que
$p\in\Qp^\times$ s'envoie sur $\mathrm{Frob}_p^{-1}$), alors on pose~:
\[ \rho(r,\chi) = (\mathrm{ind}(\omega_2^{r+1})) \otimes \chi. \]
On obtient ainsi toutes les 
repr\'esentations absolument irr\'eductibles de $\g$ de 
dimension $2$ sur $k_L$, 
et les entrelacements entre les
$\rho(r,\chi)$ sont les suivants~: 
\begin{align*} 
\rho(r,\chi) \simeq & \rho(r,\chi \mu_{-1}) \\
\rho(r,\chi) \simeq & \rho(p-1-r,\chi \omega^r) \\
\rho(r,\chi) \simeq & \rho(p-1-r,\chi \omega^r \mu_{-1}).
\end{align*} 
On en d\'eduit une bijection {\og naturelle \fg} entre les deux
classes de repr\'esentations. 

\subsection{Application aux repr\'esentations $V_{k,a_p}$}
On suppose que le Frobenius sur $D_{k,a_p}$ est semi-simple. Commen\c{c}ons par voir que la donn\'ee 
d'un r\'eseau $\Pi^0_{k,a_p,\chi}$ 
de $\Pi(V_{k,a_p,\chi})$ (i.e. d'une boule unit\'e stable
par $\G$ de ce Banach) et donc de 
$\Pi(V_{k,a_p,\chi})^* \simeq
\projlim_\psi D^\sharp(V_{k,a_p,\chi})$ 
d\'etermine un r\'eseau $T_{k,a_p,\chi}$
de $V_{k,a_p,\chi}$ stable par $\g$ tel que $(\Pi^0_{k,a_p,\chi})^* \simeq
\projlim_\psi D^\sharp(T_{k,a_p,\chi})$~:

\begin{lemm}\label{lerezo}
Si $M$ est un $\OO_L$-r\'eseau  de $\projlim_\psi
D^\sharp(V_{k,a_p,\chi})$ qui est stable par $\B$, 
alors il existe un $\OO_L$-r\'eseau
$T_{k,a_p,\chi}$ de $V_{k,a_p,\chi}$ stable par $\g$ 
et tel que $M = \projlim_\psi D^\sharp(T_{k,a_p,\chi})$. 
\end{lemm}

\begin{proof}
Soit $\pr_0 : \projlim_\psi D^\sharp(V_{k,a_p,\chi}) 
\ra D^\sharp(V_{k,a_p,\chi})$ la projection
$\{v_n\} \mapsto v_0$ et $M_0 = \pr_0(M)$. 
Par le lemme 4.57 de \cite{CL}, on a $M =
\projlim_\psi M_0$ ce qui fait que $M_0$ est un
$\OO_L[\![X]\!]$-r\'eseau stable par $\psi$ et $\Gamma$ de
$D^\sharp(V_{k,a_p,\chi})$ et donc que 
$\calO \otimes_{\OO_L[\![X]\!]} M_0$ est un 
$\calO$-r\'eseau de $D(V_{k,a_p,\chi})$. 
Par fonctorialit\'e des $(\phi,\Gamma)$-modules, il existe un 
$\OO_L$-r\'eseau $T_{k,a_p,\chi}$ de $V_{k,a_p,\chi}$ stable par $\g$ 
tel que $\calO \otimes_{\OO_L[\![X]\!]} M_0 = D(T_{k,a_p,\chi})$
et donc tel que $M = \projlim_\psi D^\sharp(T_{k,a_p,\chi})$.
\end{proof}

Par \cite[proposition 4.50]{CL}
on a un isomorphisme topologique $\B$-equivariant~: 
\[ k_L \otimes_{\OO_L} (\Pi^0_{k,a_p,\chi})^* \simeq 
\projlim_\psi D^\sharp(k_L \otimes_{\OO_L} T_{k,a_p,\chi}). \]

\begin{rema}\label{topo}
Par \cite{Br2}, $k_L \otimes_{\OO_L} 
(\Pi^0_{k,a_p,\chi})^*$ est un $k_L[\![\K]\!]$-module 
de type fini, donc muni d'une unique topologie s\'epar\'ee telle que l'action
de $k_L[\![\K]\!]$ soit continue (d\'eduite de la topologie
d'anneau noeth\'erien compact de $k_L[\![\K]\!]$). Par ailleurs, $k_L
\otimes_{\OO_L} \Pi^0_{k,a_p,\chi}$ est limite inductive d'espaces
vectoriels de dimension finie sur $k_L$ (donc finis) fixes sous
l'action de sous-groupes ouverts de plus en plus petits de $\K$. On
v\'erifie facilement que la topologie de $k_L \otimes_{\OO_L}
(\Pi^0_{k,a_p,\chi})^*$ s'identifie alors \`a la topologie de la
limite projective du dual alg\'ebrique $(k_L 
\otimes_{\OO_L} \Pi^0_{k,a_p,\chi})^*$. 
La topologie de $\projlim_\psi D^\sharp(k_L
\otimes_{\OO_L} T_{k,a_p,\chi})$ est la topologie de la limite projective,
chaque $D^\sharp(k_L \otimes_{\OO_L} T_{k,a_p,\chi})$ \'etant muni de la
topologie $X$-adique.
\end{rema}

On pose $\overline{\Pi}_{k,a_p,\chi} = 
k_L \otimes_{\OO_L} \Pi^0_{k,a_p,\chi}$. Si $2 \leq k \leq p$ et
$V_{k,a_p,\chi}$ est telle que 
$\overline{V}_{k,a_p,\chi}$ est irr\'eductible, alors
$\overline{V}_{k,a_p,\chi}$ correspond bien \`a
$\overline{\Pi}_{k,a_p,\chi}$ 
sous la bijection naturelle du paragraphe pr\'ec\'edent
et on obtient toutes les
repr\'esentations supersinguli\`eres 
de $\G$ de cette mani\`ere (voir \cite{Br2}).

\begin{lemm}\label{ddirr}
Si $U$ est une repr\'esentation irr\'eductible de $\g$ de 
dimension $2$ sur $k_L$, et si $M$ est un 
sous-$\OO_L[\![X]\!]$-module non-nul de $D^\sharp(U)$ stable par $\psi$ et
$\Gamma$, alors $M=D^\sharp(U)$. 
\end{lemm}

\begin{proof}
Commen\c{c}ons par remarquer que pour tout polyn\^ome $P \in k_L[X]$,
on a $D(U)^{P(\phi)=0} = 0$ parce que (cf. la preuve du (iii) de la 
remarque 5.5 de \cite{CL}) 
on a~: \[ D(U)^{P(\phi)=0} \subset 
(\Fpbar
\otimes_{\Fp} U)^ { \mathrm{Gal} (\Qpbar / \Qp(\mu_{p^\infty})) }
\subset \Fpbar \otimes_{\Fp} 
U^{ \mathrm{Gal} (\Qpbar / \Qp^{\mathrm{ab}}) } =  0. \]
Le lemme suit alors de la proposition 4.47 de \cite{CL} (ou plus
exactement de sa d\'emonstration, en remarquant que la d\'emonstration
n'utilise pas le fait que $\psi:M \ra M$ est surjectif). 
\end{proof}

\begin{prop}\label{galirred}
Si $U$ est une repr\'esentation irr\'eductible de $\g$ de 
dimension $2$ sur $k_L$, alors la repr\'esentation 
$\projlim_\psi D^\sharp(U)$ de $\B$ est topologiquement 
irr\'eductible. 
\end{prop}

\begin{proof}
Soit $\pr_j : \projlim_\psi D^\sharp(U) \ra D^\sharp(U)$ la projection
$\{v_n\} \mapsto v_j$. Si $M$ est un sous-espace ferm\'e et stable par
$\B$ de $\projlim_\psi D^\sharp(U)$, on note
$M_j$ l'image de $\pr_j : M \ra D^\sharp(U)$. On voit que $M_j$
est un sous-$\OO_L[\![X]\!]$-module non-nul 
de $D^\sharp(U)$ stable par $\psi$ et
$\Gamma$ ce qui fait que, par
le lemme \ref{ddirr}, $M_j=D^\sharp(U)$. On en d\'eduit que $M$ est
dense dans $\projlim_\psi D^\sharp(U)$ et donc finalement que $M =
\projlim_\psi D^\sharp(U)$ ce qui fait que 
$\projlim_\psi D^\sharp(U)$ est bien topologiquement 
irr\'eductible.
\end{proof}

\begin{prop}\label{gl2irred}
Si $\Pi$ est une repr\'esentation
supersinguli\`ere de $\G$ sur $k_L$, alors sa restriction \`a
$\B$ est irr\'eductible.
\end{prop}

\begin{proof}
Par ce qui pr\'ec\`ede et par la remarque 
\ref{topo}, le dual alg\'ebrique $\Pi^*$ de $\Pi$ avec sa topologie
profinie est topologiquement et de fa\c con $\B$-\'equivariante
isomorphe \`a $\projlim_\psi D^\sharp(\overline{V}_{k,a_p,\chi})$ pour un
$\overline{V}_{k,a_p,\chi}$ convenable avec $2\leq k\leq p$. 
Par la proposition \ref{galirred}, $\Pi^*$ est donc une
repr\'esentation topologiquement irr\'eductible de $\B$. On en
d\'eduit que $\Pi$ est irr\'eductible (un quotient strict de $\Pi$
stable par $\B$ fournirait par dualit\'e un sous-espace ferm\'e strict
de $\Pi^*$ stable par $\B$). 
\end{proof}

\begin{prop}\label{schur}
Si $U_1$, $U_2$ sont deux repr\'esentations 
irr\'eductibles de $\g$ de 
dimension $2$ sur $k_L$ et s'il existe une application
$\B$-\'equivariante continue et non-nulle $f : \projlim_\psi
D^\sharp(U_1) \ra \projlim_\psi D^\sharp(U_2)$,  
alors $U_1 \simeq U_2$. 
\end{prop}

\begin{proof}
La d\'emonstration est analogue \`a celle de la proposition 3.4.3 de
\cite{BB}. Notons comme ci-dessus $\pr_0 : \projlim_\psi D^\sharp(U)
\ra D^\sharp(U)$ la projection $\{v_n\} \mapsto v_0$.

Commen\c{c}ons par montrer que si $v=\{v_n\} \in \projlim_\psi
D^\sharp(U_1)$, alors $\pr_0 \circ f(v)$ ne d\'epend que de 
$v_0 = \pr_0(v)$. Soit $K_n$ l'ensemble des $v \in 
\projlim_\psi D^\sharp(U_1)$
dont les $n$ premiers termes sont nuls, ce qui fait 
que pour $n \geq 1$, $K_n$ est un sous-$\OO_L[\![X]\!]$-module 
ferm\'e et stable par $\psi$ et $\Gamma$ de 
$\projlim_\psi D^\sharp(U_1)$ et que $\psi(K_n)=K_{n+1}$. 
On en d\'eduit que $\pr_0 \circ f (K_n)$ est un
sous-$\OO_L[\![X]\!]$-module ferm\'e et stable par $\psi$ et $\Gamma$  
de $D^{\sharp}(U_2)$. 
Le lemme \ref{ddirr} implique alors que soit $\pr_0 \circ f (K_n) =
0$, soit $\pr_0 \circ f (K_n) = D^{\sharp}(U_2)$. Enfin, on voit que
$\psi(\pr_0 \circ f (K_n)) = \pr_0 \circ f (K_{n+1})$ et que
$\pr_0 \circ f (K_n) = 0$ si $n \gg 0$ par continuit\'e. Cela
implique que $\pr_0 \circ f (K_n) = 0$ pour tout $n \geq 1$ et donc
que si $v_0 = 0$, alors $\pr_0 \circ f(v) = 0$. 

Pour tout $w \in D^{\sharp}(U_1)$, soit $\widetilde{w}$ un \'el\'ement de 
$\projlim_{\psi} D^{\sharp}(U_1)$ tel que $\widetilde{w}_0 = w$. Les
calculs pr\'ec\'edents montrent que l'application 
$h:  D^{\sharp}(U_1) \rightarrow D^{\sharp}(U_2)$ donn\'ee par
$h(w) = \pr_0 \circ f (\widetilde{w})$ est bien d\'efinie, et
qu'elle est $\OO_L[\![X]\!]$-lin\'eaire et commute \`a $\psi$ et \`a l'action
de $\Gamma$. Par les propositions 4.7 et 4.55 de 
\cite{CL}, elle s'\'etend en une application de 
$(\phi,\Gamma)$-modules $h: D(U_1) \ra D(U_2)$ et par fonctorialit\'e,
on en d\'eduit qu'il existe une application non-nulle et
$\g$-\'equivariante de $U_1$ dans $U_2$, ce qui fait que par le lemme
de Schur, on a $U_1 \simeq U_2$.
\end{proof}

\begin{coro}\label{fin3}
Si $2p \geq k \geq p+2$ et $0 < \val(a_p) < 1$, alors 
$\overline{V}_{k,a_p} = \mathrm{ind} (\omega_2^{k-p})$. 
\end{coro}

\begin{proof}
Notons d'abord que, sous les conditions de l'\'enonc\'e, 
le Frobenius sur $D_{k,a_p}$ est toujours semi-simple (v\'erification facile).
La proposition 5.8 de \cite{Br1} montre que~: 
\[ \overline{\Pi}_{k,a_p} \simeq \pi(2p-k,0,\omega^{k-1-p}) =
\pi(k-p-1,0,1) \]
et la proposition \ref{gl2irred} (ou plut\^ot 
sa preuve) entra\^\i ne que la restriction de 
$\overline{\Pi}_{k,a_p} ^*\simeq \projlim_\psi
D^\sharp(\overline{V}_{k,a_p})$ \`a $\B$ est topologiquement
irr\'eductible. On en d\'eduit que $\overline{V}_{k,a_p}$ est
elle-m\^eme irr\'eductible (s'il existait une sous-repr\'esentation
stricte $U$ de $\overline{V}_{k,a_p}$, on en d\'eduirait l'existence d'un
sous-espace ferm\'e $\projlim_\psi D^\sharp(U)$ 
de $\projlim_\psi D^\sharp(\overline{V}_{k,a_p})$ stable sous $\B$). 

On a d'autre part un isomorphisme $\projlim_\psi
D^\sharp(\rho(k-p-1,1)) \simeq \pi(k-p-1,0,1)^*$ (par la proposition 6.2 de
\cite{Br2} par exemple) et donc un isomorphisme 
topologique et $\B$-\'equivariant $\projlim_\psi
D^\sharp(\rho(k-p-1,1)) \simeq \projlim_\psi
D^\sharp(\overline{V}_{k,a_p})$ ce qui fait, par la proposition
\ref{schur}, que l'on a bien $\overline{V}_{k,a_p} \simeq 
\rho(k-p-1,1) = \mathrm{ind} (\omega_2^{k-p})$.
\end{proof}

\begin{rema}
Ce corollaire est un cas particulier de la conjecture 
selon laquelle lorsque l'on r\'eduit modulo $p$,
la correspondance $V_{k,a_p,\chi} \leftrightarrow
\Pi_{k,a_p,\chi}$ est compatible avec la correspondance modulo $p$
rappel\'ee plus haut (et \'etendue aux cas r\'eductibles).
\end{rema}


\begin{thebibliography}{Bre03b}
 
\bibitem[Ben00]{Bn1} \textsc{Benois, D :} 
\textit{On Iwasawa theory of crystalline representations},
Duke Math. J. 104 (2000), no. 2, 211--267. 

\bibitem[Ber02]{Be1} \textsc{Berger, L :} 
\textit{Limites de repr\'esentations cristallines},   
Compos. Math.  140  (2004),  no. 6, 1473--1498.

\bibitem[BB04]{BB} \textsc{Berger, L; Breuil, C :}
\textit{Repr\'esentations cristallines irr\'eductibles de ${\rm
    GL}_2(\Qp)$}, 
pr\'epublication, octobre 2004.

\bibitem[BLZ04]{BLZ} \textsc{Berger, L; Li, H; Zhu, H :}
\textit{Construction of some families of $2$-dimensional 
crystalline representations},  
Math. Ann. 329 (2004), no. 2, 365--377.

\bibitem[Bre03a]{Br1} \textsc{Breuil, C :} 
\textit{Sur quelques repr\'esentations modulaires 
et $p$-adiques de ${\rm GL}_2(\Qp)$ I}, 
Compositio Math.  138  (2003),  no. 2, 165--188.

\bibitem[Bre03b]{Br2} \textsc{Breuil, C :} 
\textit{Sur quelques repr\'esentations modulaires 
et $p$-adiques de ${\rm GL}_2(\Qp)$ II}, 
J. Institut Math. Jussieu 2, 2003, 23--58.

\bibitem[CC99]{CC} \textsc{Cherbonnier, F ; Colmez, P :} 
\textit{Th\'eorie d'Iwasawa des repr\'esentations $p$-adiques d'un corps local},  
J. Amer. Math. Soc. 12 (1999), no. 1, 241--268.

\bibitem[Col04]{CL} \textsc{Colmez, P :} 
\textit{Une correspondance de Langlands locale $p$-adique pour les
  repr\'esen\-tations semi-stables de dimension $2$},
pr\'epublication, 2004.

\bibitem[CF00]{CF} \textsc{Colmez, P;  Fontaine, J-M :} 
\textit{Construction des repr\'esentations $p$-adiques semi-stables}, 
Inv. Math. 140, 2000, 1--43.

\bibitem[Fon90]{F90} \textsc{Fontaine, J-M :} 
\textit{Repr{\'e}sentations $p$-adiques des corps locaux I}, 
The Grothendieck Festschrift,
Vol. II, 249--309, Progr. Math. 87, Birkh{\"a}user Boston, Boston,
MA 1990.

\bibitem[FL82]{FL82}
\textsc{Fontaine, J-M; Laffaille G :} 
\textit{Construction de repr\'esentations $p$-adiques},  
Ann. Sci. \'Ecole Norm. Sup. (4) 15 (1982), no. 4, 547--608 (1983). 

\bibitem[Wa96]{W96} 
\textsc{Wach, N :} 
\textit{Repr{\'e}sentations $p$-adiques potentiellement cristallines}, 
Bull. Soc. Math. France 124, 1996, 375--400.
\end{thebibliography}
\end{document}